\documentclass[a4paper,11pt,oneside]{amsart}
\usepackage[T1]{fontenc}
\usepackage{microtype}
\usepackage[hmargin=3cm,vmargin=2.5cm]{geometry}
\usepackage{amsmath,amssymb,amsthm,mathtools}
\usepackage{mathrsfs}
\allowdisplaybreaks[2]
\usepackage{graphicx}
\graphicspath{{figures/}}
\usepackage[dvipsnames]{xcolor}
\usepackage[cmtip,all]{xy}
\usepackage{tikz}
\usetikzlibrary{arrows.meta,calc,positioning}
\usepackage{float}
\usepackage{booktabs}
\usepackage{enumitem}
\usepackage{appendix}
\tikzset{
  vtx/.style={
    circle,
    draw,
    line width=0.6pt,
    fill=white,
    minimum size=3mm,
    inner sep=0pt
  },
  edg/.style={line width=0.6pt}
}
\usepackage{csquotes}
\usepackage[
  backend=biber,
  style=alphabetic,
  sorting=nyt,
  giveninits=true,
  maxbibnames=4,
  maxcitenames=3,
  doi=false,
  isbn=false,
  url=false
]{biblatex}
\usepackage[
    colorlinks=true,
    linkcolor=MidnightBlue,
    citecolor=ForestGreen,
    urlcolor=BrickRed]{hyperref}
\theoremstyle{plain}
\newtheorem{theorem}{Theorem}[section]

\newtheorem{proposition}[theorem]{Proposition}
\newtheorem{corollary}[theorem]{Corollary}

\newtheorem{fact}[theorem]{Fact}

\theoremstyle{definition}
\newtheorem{definition}[theorem]{Definition}

\theoremstyle{remark}
\newtheorem{remark}[theorem]{Remark}

\numberwithin{equation}{section}
\newcommand{\Q}{\mathbb{Q}}
\newcommand{\R}{\mathbb{R}}
\newcommand{\Z}{\mathbb{Z}}
\newcommand{\A}{\mathbb{A}}

\renewcommand{\P}{\mathbb{P}}

\newcommand{\Acal}{\mathcal{A}}
\newcommand{\Ccal}{\mathcal{C}}

\newcommand{\vb}{\mathbf{v}}
\newcommand{\xb}{\mathbf{x}}
\newcommand{\Am}{\mathscr A_{\boldsymbol m}}
\newcommand{\Bm}{\boldsymbol m}
\newcommand{\Bc}{\boldsymbol c}

\DeclareMathOperator{\GL}{GL}
\DeclareMathOperator{\SL}{SL}

\DeclareMathOperator{\Mat}{Mat}

\DeclareMathOperator{\rank}{rank}
\DeclareMathOperator{\diag}{diag}
\DeclareMathOperator{\tr}{tr}

\DeclareMathOperator{\rk}{rk}

\newcommand{\Id}{\mathrm{Id}}

\title[Infinitude of integral Hitchin representations of Coxeter polygons]{Infinitude of integral Hitchin representations of Coxeter polygons}

\author[S. Kim]{Seong youn Kim}
\address{Department of Mathematical Sciences, Seoul National University, Seoul 08826, South Korea}
\email{s264575@snu.ac.kr}
\author[S. Ko]{Sunghwan Ko}
\address{Department of Mathematical Sciences, Seoul National University, Seoul 08826, South Korea}
\email{sym2346@snu.ac.kr}

\subjclass[2020]{22E40, 20F55, 13F60}
\keywords{Integral representations, Hitchin component, Coxeter polygons, cluster algebras, generalized Markoff equations}
\date{\today}

\begin{document}

\begin{abstract}
We prove that there are infinitely many Hitchin $\SL^\pm(3, \Z)$-representations of a Coxeter $n$-gon($n \geq 9$) group whose boundary parameters are fixed. We construct infinitely many Hitchin $\SL^\pm(3, \Z)$-representations of a specific Coxeter $6$-gon or $8$-gon corresponding to positive integral points of a generalized Markoff equation.
\end{abstract}

\maketitle


\section{Introduction}\label{sec:introduction}

\subsection*{Main results}
Let $P$ be a Coxeter $n$-gon with cyclically ordered labels ${\Bm}=(m_1,\dots,m_n)$ and let $W_{\Bm}$ be its Coxeter group. A representation
\[
\rho:W_{\Bm} \longrightarrow \SL^\pm(3,\R)
\]
is called \emph{Vinberg} if it is induced by the side reflections of a properly convex projective Coxeter polygon realizing $P$. If $m_i$ is finite for every $i$ and $\sum_{i=1}^{n}1/m_i<n-2$, the space of Vinberg representations up to conjugacy coincides with the $\SL^\pm(3,\R)$-Hitchin component associated with $P$ \cite{ALS23}. If some of $m_i$s are infinite, then $\rho$ is called \textit{Hitchin} if its double representation into $\SL(3, \R)$ is cusped Hitchin; cf.\cite{canary2022cusped}. Direct computation shows that any Vinberg representation is Hitchin. A Vinberg representation is called \emph{integral} if it is conjugate into $\SL^\pm(3,\Z)$. We ask whether integral Vinberg representations with fixed boundary parameters form finitely or infinitely many conjugacy classes.

Integral discrete faithful representations in geometric connected components of character varieties have been studied through both constructions of infinite families and finiteness theorems. In \cite{LRT2011}, they constructed infinitely many pairwise non-conjugate integral representations of triangle group $\Delta(3,4,4)$ in its Hitchin component, with Zariski-dense images in $\SL(3,\R)$. In contrast, it is well-known that every Coxeter $3$-gon group admits only finitely many integral representations, up to conjugacy, when its boundary parameters are fixed; cf.\cite{CC15}. More recently, Ko proved finiteness of integral representations associated with irreducible, large, $2$-perfect Coxeter truncation polytopes of dimension at least three \cite{Ko26}. 

Note that every Coxeter polygon is combinatorially a truncation polytope, since it can be obtained from a triangle by successive vertex truncations. Our main first theorem shows that the above finiteness result does not extend to dimension two in general.

\begin{theorem}[Theorem~\ref{infinite Cartan}]\label{mainthm1}
    Let $P$ be a Coxeter $n$-gon with ${\Bm}=(m_1,\dots,m_n)$ and let $W_{\Bm}$ be its Coxeter group. Suppose $n\ge 9$ and $m_i\in\{2,3,4,6,\infty\}$ for every $i$. For every choice of integers $c_i\ge4$ at the vertices labeled by $\infty$, $W_{\Bm}$ admits infinitely many non-conjugate integral Vinberg representations in Hitchin component.
\end{theorem}

A direct application of Theorem~\ref{mainthm1} is the construction of infinite families of thin reflection groups. Recall that a finitely generated subgroup of $\SL^\pm(3,\Z)$ is called \emph{thin} if it is Zariski-dense in $\SL^\pm(3,\R)$ and has infinite index in $\SL^\pm(3,\Z)$.

Indeed, suppose that all labels in $\Bm$ are finite, then $P$ is compact. Theorem~\ref{mainthm1} provides infinitely many  integral representations in Hitchin component associated with $P$. Takeuchi proved that arithmetic Fuchsian groups of a fixed signature form finitely many $\SL(2,\R)$-conjugacy classes \cite{Tak83}. Since integral Fuchsian representations in Hitchin component have arithmetic images, and it implies that only finitely many integral conjugacy classes lie in the Fuchsian locus. 

In \cite{AudibertDoubaLeeMarquis25_ZariskiClosures}, they proved that a representation in Hitchin component of $P$ has Zariski-dense image in $\SL^\pm(3,\R)$ if and only if it is not contained in Fuchsian locus. Consequently, all but finitely many of integral conjugacy classes above have Zariski-dense images and hence they are thin.

The main object in the proof is the Grassmannian cluster algebra $\Acal_{3,n}$ \cite{scott2006grassmannians}. In \cite{MOT12}, they established a correspondence between closed $2$-friezes and normalized vector polygons, relating their geometry to cluster algebras. Morier-Genoud subsequently proved that there are infinitely many positive integral closed $2$-frieze of width $n-4$ for $n\ge9$ \cite{MG12}.

Our construction extends this geometric interpretation by constructing projective reflections in the sides of the associated dual polygons. Starting from normalized polygon data, we use the unimodularity of consecutive triples to reconstruct integral reflection vector and hence we obtain integral reflections realizing the Coxeter polygon $P$. These reflections define integral Vinberg representations of $W_{\Bm}$. For fixed Coxeter data, integrality ensures that the above reconstruction has finite fibers after passing to equivalence classes under the positive diagonal action. The infinitude of positive integral points $\Acal_{3,n}$ gives infinitely many pairwise non-conjugate integral Vinberg representations.

For $n<9$, on the other hand, the cluster algebra $\Acal_{3,n}$ is of finite type, so the preceding infinitude argument does not apply. Our second main theorem establishes both finiteness and infinitude for smaller polygons.

\begin{theorem}\label{mainthm2}
    Let $P$ be a Coxeter $n$-gon with ${\Bm}=(m_1,\dots,m_n)$ and let $W_{\Bm}$ be its Coxeter group.
    \begin{enumerate}
        \item If $n=5$, then $W_{\Bm}$ admits only finitely many integral Vinberg representations up to conjugacy with fixed boundary parameters.
        \item If $n\in\{6,8\}$ and $m_i=3$ for every $i$, then $W_{\Bm}$ admits infinitely many pairwise non-conjugate integral Vinberg representations.
    \end{enumerate}
\end{theorem}

For hexagons and octagons, we construct families of representations in Hitchin components. These families are parametrized by positive real solutions of generalized Markoff equations. Such equations classically describe relative $\SL(2,\mathbb C)$-character varieties of the four-holed sphere \cite{CL09,whang2020nonlinear}.

We construct infinite orbits of positive integral points under the group generated by the three Vieta involutions. These points determine pairwise non-conjugate integral Vinberg representations. Moreover, these orbits are Zariski-dense in the surfaces defined by the corresponding generalized Markoff equations \cite{CL09}. In particular, these integral representations are not contained in any algebraic curve. This distinguishes our construction from the explicit one-parameter constructions of integral Hitchin representations in \cite{LRT2011,DownsThistlethwaite}, where integral representations lie on algebraic curves in the corresponding Hitchin components.

\subsection*{Remarks on integral points} 
From a Diophantine viewpoint, our problem can be reformulated as the study of integral points in a (semi)algebraic locus of an algebraic variety. It is natural to ask whether they admit a notion of finite generation: can every integral solution be obtained from finitely many initial solutions by a specified collection of natural transformations?

If all of $m_i$s were finite, the outer automorphism group $\mathrm{Out}(W_{\Bm})$ is finite. Consequently, Theorem \ref{mainthm1} yields infinitely many orbits of integral Vinberg representations under this natural action. In this sense, our result provides an analogue for Coxeter polygon groups of Audibert’s constructions of infinitely many mapping class group orbits of Hitchin representations in arithmetic lattices \cite{audibert2025zariski}.

Our cluster algebra construction suggests a broader framework for this question. Cluster mutations are birational coordinate transformations on the auxiliary Grassmannian parameter spaces; they need not induce automorphisms over $\Z$ of the character variety or bijections of the integral locus in the Hitchin component. Nevertheless, the Laurent phenomenon allows mutations to be used in constructing positive integral points. This suggests asking whether a suitable construction using mutations can produce all integral representations from a finite collection of initial data.

A necessary condition for such a description is that only finitely many tuples of frozen parameter values occur. The $6$-gon or $8$-gon family in Theorem \ref{mainthm2} shows that this finiteness can fail: its normalized frozen parameters take infinitely many distinct values, even though the associated Grassmannian cluster algebra $\Acal_{3,6}$ is of finite type. These phenomena suggest that the arithmetic of higher-rank local systems can differ substantially from that in the $\SL_2$-case. 

\subsection*{Organization of the paper}

In Section \ref{sec:prelim}, we recall the description of Vinberg representations in terms of Cartan matrices and cyclic products. We also review the Grassmannian cluster algebra and the realization of a Coxeter polygon. In Section \ref{sec:proof-main}, we prove Theorem \ref{mainthm1} by showing that the reconstruction map has finite fibers after passing to conjugacy classes. Section \ref{sec:proof-main2} is devoted to the proof of Theorem \ref{mainthm2} with establishing the finiteness of the integral points for pentagon cases.

\subsection*{Acknowledgement}
The authors thank Gye-Seon Lee and Junho Peter Whang. The first author was partially supported by Samsung Science and Technology Foundation under Project Number SSTF-BA2201-03. The second author was supported by the National Research Foundation of Korea(NRF) grant funded by the Korea government(MSIT) (RS-2023-00252171).

\subsection*{AI disclosure}
The authors used ChatGPT in developing the ideas presented in this paper. Starting from an idea suggested by the tool for the right-angled case, we obtained our first main result, Theorem~1.1, through further discussions with the tool. The discovery of the results for $n<9$ presented in this paper was first suggested by Gye-Seon Lee's discussion with the tool, as an infinite family of integral points. The authors then developed these ideas 
using generalized Markoff equations.

\section{Preliminaries}\label{sec:prelim}

\subsection{Vinberg representations and cyclic products}

Let $Q$ be a labeled polygon with $n\ge3$ sides $e_1,\dots,e_n$ indexed cyclically along $\partial Q$. Denote by $m_i\in\{2,3,\dots,\}\cup\{\infty\}$ the label of the vertex $e_i\cap e_{i+1}$, and write $\boldsymbol m=(m_1,\dots,m_n)$. Let $S=\{s_1,\dots,s_n\}$ correspond to the edges of $Q$, and let $(W_{\boldsymbol m},S)$ be the associated Coxeter system, where

\[
W_{\Bm}=\left\langle s_1,\dots,s_n\;\middle|\;s_i^2=1,\;(s_is_{i+1})^{m_i}=1\text{ whenever } m_i<\infty\right\rangle.
\]

Let $P\subset \P(\R^3)$ be a properly convex projective $n$-gon with sides $F_1,\dots,F_n$ corresponding to $e_1,\dots,e_n$ in cyclic order, up to reversal. Choose a cone lift $\widetilde P$ of $P$ and linear functional $\alpha_i\in(\R^3)^*$ such that
\[
\widetilde P=\{x\in\R^3:\alpha_i(x)\le0\text{ for every } i\},\qquad F_i=P\cap\P(\ker\alpha_i).
\]
For vectors $b_i\in\R^3$ satisfying $\alpha_i(b_i)=2$, set $\sigma_i(x)=x-\alpha_i(x)b_i$. Each $\sigma_i$ is a linear reflection with fixed hyperplane $\ker\alpha_i$.

\begin{definition}\label{def:coxeter-polygon}
    The tuple $(P,(\sigma_i)_{1\le i\le n})$ is a \emph{Coxeter polygon realizing $Q$} if, for all distinct $i,j$, $\alpha_i(b_j)\le 0$ and $\alpha_i(b_j)=0$ if and only if $\alpha_j(b_i)=0$, and if
    \[
    \alpha_i(b_{i+1})\alpha_{i+1}(b_i)=4\cos^2\!\left(\frac{\pi}{m_i}\right)\qquad\text{whenever }m_i<\infty,
    \]
    while $\alpha_i(b_j)\alpha_j(b_i)\ge4$ whenever $i,j$ are distinct and nonconsecutive, or $\{i,j\}=\{\ell,\ell+1\}$ for some $\ell$ with $m_\ell=\infty$.
\end{definition}

For the Coxeter polygon $P$ realizing $Q$, the matrix $A=(\alpha_i(b_j))_{1\le i,j\le n}$ is called a \emph{Cartan matrix} of $(P,(\sigma_i)_{1\le i\le n})$. A matrix $A\in M_n(\R)$ is \emph{compatible with $Q$} if it is a Cartan matrix of some Coxeter polygon realizing $Q$.

Such a matrix is \emph{irreducible} if it cannot be transformed into a nontrivial block diagonal matrix by simultaneously rearranging its rows and columns. If $A$ is irreducible, write $A=2\Id-B$ with $B$ nonnegative. We say that $A$ is of \emph{negative type} if the spectral radius of $B$ is greater than $2$.

Let $\Am$ denote the space of irreducible matrices compatible with $Q$ that are of negative type and rank three. We record the following elementary observation.
\begin{fact}\label{polygon-cartan}
Suppose that $n\geq5$. Then every matrix compatible with $Q$ is irreducible, of negative type, and of rank three. In particular, $\Am$ consists of all matrices compatible with $Q$.
\end{fact}
 
Set $I_\infty=\{i\in\Z/n\Z : m_i=\infty\}$, and fix a collection of real numbers $\Bc=(c_i)_{i\in I_\infty}$ with $c_i\ge4$. We define the subspace 
\[
\Am(\Bc)=\left\{A\in\Am\;\middle|\;a_{i,i+1}a_{i+1,i}=c_i\text{ for every } i\in I_\infty\right\}.
\]
Thus, the products associated with consecutive indices are determined by the Coxeter exponent $m_i$ when $m_i<\infty$, and by the parameter $c_i$ when $m_i=\infty$.

Two matrices $A,A'\in\Am$ are said to be \emph{equivalent} if
\[
A'=DAD^{-1}
\]
for some positive diagonal matrix $D$. Since each product $a_{ij}a_{ji}$ is invariant under positive diagonal conjugation, this action preserves $\Am(\Bc)$. We write
\[
\Ccal(\Bm;\Bc):=\Am(\Bc)/\!\sim.
\]

\begin{remark}
    If all $m_i$ are finite, then $I_\infty=\varnothing$, so no additional parameters $\Bc$ are needed. Moreover, $\Bm$ is realized by a compact hyperbolic Coxeter $n$-gon $P_0$ when $\sum_{i=1}^n 1/m_i<n-2$.
\end{remark}

The following is a form of the Tits--Vinberg theorem; see \cite{Bou02} and \cite[Theorem~2]{Vinberg71a}.

\begin{theorem}[Tits--Vinberg]\label{thm:tits--vinberg}
    Let $(P,(\sigma_i)_{1\le i\le n})$ be a Coxeter polygon realizing $Q$. The assignment $s_i\mapsto\sigma_i$ extends to a discrete and faithful representation
    \[
    \rho: W_{\Bm}\longrightarrow\SL^\pm(3,\R).
    \]
\end{theorem}

For $A\in\Am$, Choose a Coxeter polygon realizing $Q$ with $(\alpha_i,b_i)_{1\le i \le n}$ such that $A=(\alpha_i(b_j))_{1\le i,j\le n}$. Denote by $\rho_A$ the representation given by Theorem~\ref{thm:tits--vinberg}, so that 
\[
\rho_A(s_i)(x)=x-\alpha_i(x)b_i.
\]
The rank condition implies $\operatorname{span}\{b_1,\dots,b_n\}=\R^3$ and $\bigcap_i\ker\alpha_i=\{0\}$. Since $A$ is irreducible, $\rho_A$ is also irreducible by \cite[Proposition~3.23]{DGKLM25}.

For $g\in\mathrm{GL}(3,\R)$, the change of coordinates $(\alpha_i,b_i)\mapsto(\alpha_i\circ g^{-1},gb_i)$ preserves $A$ and conjugates $\rho_A$ by $g$. For $\lambda_1,\dots,\lambda_n>0$, rescaling $(\alpha_i,b_i)\mapsto(\lambda_i\alpha_i,\lambda_i^{-1}b_i)$ preserves $\rho_A$, and replaces $A$ by $DAD^{-1}$, where $D=\diag(\lambda_1,\dots,\lambda_n)$. Uniqueness of rank factorization up to change of basis implies $[A]\mapsto[\rho_A]$ is well defined.

A representation $\rho\colon W_{\boldsymbol m}\longrightarrow\SL^\pm(3,\R)$ is called a \emph{Vinberg representation of type $\boldsymbol m$} if it is conjugate to $\rho_A$ for some $A\in\Am$. Let $\chi\bigl(W_{\boldsymbol m},\SL^\pm(3,\R)\bigr)$ denote the space of conjugacy classes of semisimple representations. Define
\[
\mathscr X_{\boldsymbol m}^{\mathrm{Vin}} :=\left\{ [\rho_A]\in\chi\bigl(W_{\boldsymbol m},\SL^\pm(3,\R)\bigr) \;\middle|\; A\in\Am\right\}.
\]
Since every $\rho_A$ is irreducible, it is semisimple. By \cite[Fact~3.29]{DGKLM25} and the previous observation, the assignment
\[
\Phi\colon\Am/\!\sim\longrightarrow\mathscr X_{\boldsymbol m}^{\mathrm{Vin}},\qquad[A]\longmapsto[\rho_A],
\]
is a homeomorphism.

For $\Bc=(c_i)_{i\in I_\infty}$ as above, define
\[
\mathscr X_{\boldsymbol m}^{\mathrm{Vin}}(\Bc):= \left\{ [\rho]\in\mathscr X_{\boldsymbol m}^{\mathrm{Vin}} \;\middle|\; \operatorname{tr}\bigl(\rho(s_i s_{i+1})\bigr)=c_i-1 \text{ for every }i\in I_\infty \right\}.
\]
The identity $\operatorname{tr}\bigl(\rho_A(s_is_{i+1})\bigr) =a_{i,i+1}a_{i+1,i}-1$ implies that $\Phi$ restricts to a homeomorphism
\[
\Phi_{\Bc}\colon \Ccal(\Bm;\Bc)\longrightarrow \mathscr X_{\boldsymbol m}^{\mathrm{Vin}}(\Bc), \qquad [A]\longmapsto[\rho_A].
\]

We now introduce the cyclic-product invariants. For $\ell\ge2$, an \emph{oriented circuit} in $S$ is a tuple $C=(s_{i_1},\dots,s_{i_\ell})$ of pairwise distinct elements of $S$. Its corresponding \emph{cyclic product} is
\[
C(A):= a_{i_1i_2}a_{i_2i_3}\cdots a_{i_\ell i_1}.
\]

Since each diagonal factor appearing in $DAD^{-1}$ cancels in the product, $C(A)$ depends only on the equivalence class $[A]$.

\begin{theorem}[{\cite[Prop.~16]{Vinberg71a}}]
\label{thm:vinberg-cyclic-products}
Let $A,A'\in\Am$. Then $A$ and $A'$ are equivalent if and only if
\[
C(A)=C(A')
\]
for every oriented circuit $C$ in $S$.
\end{theorem}

We call $\rho_A$ \emph{integral} if it preserves a lattice in $\R^3$, or equivalently, if it is conjugate to a representation with image in $\SL^\pm(3,\Z)$. In the present setting, \cite[Lemma~6.2]{AudibertDoubaLeeMarquis25_ZariskiClosures} takes the following form.

\begin{proposition}[Integrality criterion]\label{Integrality criterion}
Let $A \in \Am$. Then $\rho_A$ is integral if and only if
\[
C(A)\in\Z
\]
for every oriented circuit $C$ in $S$.
\end{proposition}

Consequently, in order to construct integral Vinberg representations of $W_{\boldsymbol m}$, it is enough to construct a matrix $A\in\Am$ whose cyclic products are all integers.

We now slightly strengthen the statement. We call a vector(resp. covector) $v \in \Z^3 \subseteq \R^3$(resp. $(\Z^3)^* \subseteq (\R^3)^*$) is \textit{primitive} if its entries are coprime.

\begin{proposition}\label{IntegralCartan}
    Let $A \in \Am$. Then $\rho_A$ is integral if and only if there exists $A' \sim A$ such that
    \[
    A' \in \Mat_{n \times n}(\Z).
    \]
    Moreover, for a fixed $A$ the number of such $A'$ is finite, only depending on $\Bc$.
\end{proposition}

\begin{proof}
    It suffices to check the only if direction. Suppose $\rho_A(s_i) \in \SL^\pm(3, \Z)$ for any $s_i$. Then one can choose the covector $\alpha_i \in (\Z^3)^*$ to be primitive. Since $\alpha_i$ is primitive, there exists $v_i \in \Z^3$ such that $\alpha_i(v_i) = 1$. Then
    \[
    b_j = \alpha_j(v_j)b_j = (I - \sigma_j)v_j \in \Z^3
    \]
    and $\alpha_i(b_j) \in \Z$ for any $i, j$. For a fixed $A$, the $2$-cyclic product
    \[
    \alpha_i(b_j) \alpha_j(b_i) \in \Z
    \]
    is fixed. Thus any entry of $A'$ should divide a fixed constant only depending on $\Bc$, whose possibility is finite.
    
\end{proof}

By the above proposition, we can write $A \in \Am(\Z)$ if $\rho_A$ is integral. Denote 
\[
\Ccal(\Bm ; \Bc)(\Z) := \{ A \in \Am : A \sim A' \text{ for some } A' \in \Am(\Z) \} / \sim.
\]

\subsection{Affine Grassmannian cluster algebra}

Let $K$ be a field. For the simplicity, we will only concern the case $K = \R$ but can be generalized for arbitrary (ordered) field $K$ (of characteristic $0$). Let
\[
M = ( \vb_1 , \vb_2 , \cdots , \vb_n ) \in \Mat_{3 \times n}(\R)
\]
where $\vb_i \in \R^3$. Consider each $\vb_i$ as a vector-valued regular function on $\Mat_{3 \times n}(k) \simeq k^{3n}$. For a $3$-tuple $I = (i, j, k) \in [n]^3$, the regular function
\[
\Delta_I(M) = \det(\vb_i, \vb_j, \vb_k)
\]
is called the \textit{$3$-Pl\"ucker coordinate function}. From now on, we may abuse $\Delta_I(M) = \Delta_I$ if the target matrix $M$ has no ambiguity. There is a natural \textit{$3$-Pl\"ucker relation}
\[
\Delta_{i_1, i_2, i_3} \Delta_{j_1, j_2, j_3} = \sum_{r=1}^3 \Delta_{i_1, i_2, j_r}
\Delta_{j_1, \ldots, j_{r-1}, i_3, j_{r+1}, \ldots, j_3}.
\]
The \textit{$3$-Pl\"ucker algebra} is given by
\[
\Acal_{3,n} := \R\left[ \Delta_I : I \in [n]^3 \right].
\]
Write $\Delta_i = \Delta_{i+1, i+2, i+3}$ where we identify $\Z / n\Z$ with $[n]$ and use the cyclic order. We will call the tuple $(i, j, k)$ to be \textit{positive} if it admits a lift to $\Z / n\Z$ with $i < j < k < i+n$.
 
\begin{theorem}[{\cite[Theorem~3]{scott2006grassmannians}}]
    $\Acal_{3,n}$ is a cluster algebra whose (initial) cluster variables consisting of $(2n-8)$ many $\Delta_I$s with frozen variables $\Delta_i$.
\end{theorem}

The choice of the $(2n-8)$ variables is obtained by the \textit{Postnikov diagram} \cite{scott2006grassmannians}. By the Laurent phenomenon of the cluster algebra, we have 
\[
\Acal_{3, n} \subseteq \R[\Delta_0, \ldots, \Delta_{n-1}][\Delta_{I_s}^{\pm}]_{s=1}^{2n-8}.
\]

Consider the basic open affine chart $D( \Delta_0 )$ of $\Acal_{3,n}$ defined by $\Delta_0 \neq 0$. Gaussian elimination tells 
\[
(\vb_1, \vb_2, \vb_3)^{-1}M=(I_3\mid X), \quad
X \in \Mat_{3 \times (n-3)}(\R)
\]
and then $D(\Delta_0) \simeq \A_\R^{3(n-3)}$. We will abuse the homogenization $\Delta_I / \Delta_0$ by $\Delta_I$, giving the cluster algebra 
\[
\overline{\Acal_{3,n}} = \Acal_{3,n} / (\Delta_0 - 1) \subseteq \R\left[ \Delta_1, \ldots, \Delta_{n-1} \right]\left [\Delta_{I_s} \right]_{s=1}^{2n-8} .
\]
Let $\xb = (x_1, \ldots, x_{2n-8}; \Delta_1, \ldots, \Delta_{n-1})$ be a cluster of $\overline{\Acal_{3,n}}$. Denote the nonvanishing locus of its all cluster variables as
\[
D(\xb) := D( x_1 \cdots x_{2n-8}\Delta_1 \cdots \Delta_{n-1}).
\]

\begin{theorem}[{\cite[Corollary~4]{scott2006grassmannians}}]\label{ScottTorus}
    For any cluster $\xb$ of $\Acal_{3,n}$, we have an isomorphism of a variety over $\R$
    \[
    \Phi_{\xb} : D(\xb) \xrightarrow{\sim} (\mathbb{G}_m)_\R^{3(n-3)}.
    \]
    Consider the initial cluster $\mathbf{\Delta} = (\Delta_{I_1}, \ldots, \Delta_{I_{2n-8}} ;  \Delta_0, \ldots, \Delta_{n-1})$, we have a bijection
    \[
    D^+( \mathbf{\Delta} ) = \left\{ M \in \Mat_{3 \times (n-3)}(\R) : \Delta_I(I_3 \mid M) > 0, \forall I \text{ positive $3$-tuple} \right\} \simeq \R_{>0}^{3(n-3)}.
    \]
\end{theorem}

For the bijection, observe that the matrix $(I_3 \mid M) = (m_{ij})$ is given by 
\[
\mathbf{m}_j = (m_{1j}, m_{2j}, m_{3j}) = (\Delta_{(2, 3, j)}, -\Delta_{(1, 3, j)}, \Delta_{(1, 2, j)})(I_3 \mid M).
\]

\subsection{Arithmeticity}
For $A \in \Am$, write the rank factorization as
\[
R=(\alpha_1, \cdots, \alpha_n),
\quad U=(b_1, \cdots, b_n),\quad A=R^tU, \quad \rk R = 3 = \rk U.
\]
$R, U$ considered as an ordered tuple in $(\R^3)^n$ is called a \textit{realization} of $A$. We may write $DR \sim R$ for any positive diagonal matrix $D$ since $DAD^{-1} = (DR)^t(DU)^{-1}$. In particular 

\begin{proposition}
    Let $A \in \Am(\Z)$. Then we can take a realization $R, U \in (\Z^3)^n$. Moreover, $R, U$ is a point of $\Q_{>0}^{3(n-3)}$ via the bijection in Theorem \ref{ScottTorus}.
\end{proposition}

\begin{proof}
    That $R \in (\Z^3)^n$ is shown in the proof of Proposition \ref{IntegralCartan}. Since $P_A$ is properly convex, its dual
    \[
    (P_A)^* = \mathbb{P}(\mathcal{CH}(R))
    \]
    is also properly convex where $\mathcal{CH}$ denotes the convex hull. By replacing $\alpha_i$ by $- \alpha_i$ if we need, we may assume $\Delta_I(R) > 0$ for any positive $3$-tuple $I$ and $R \in D^+(\mathbf{\Delta})$. Note the bijection is obtained by the homogenized variable $\Delta_I(R) / \Delta_0(R)$, which should be a rational.
    
    Since $A^t \in \Am(\Z)$, the statement also symmetrically holds for $U$.

\end{proof}

    Let 
    \begin{align*}
        \mathscr{G}_{\Bm}(\Bc)(\Z) := \{R \in \Mat_{3 \times n}(\R) : &\exists A \in \Am(\Bc)(\Z) \text{ realized by } \\
        &R \sim (I_3 \mid M) \text{ for some } M \in D^+(\mathbf{\Delta})(\Q) \}.
    \end{align*}

\begin{theorem}\label{IntegralFactorizationII}
    Let $R \in \Mat_{3 \times n}(\R)$. Then there exists a finite-time algorithm testing whether $R \in \mathscr{G}_{\Bm}(\Bc)(\Z)$ or not.
\end{theorem}

\begin{proof}
    We may assume $R = (I_3 \mid M)$ and there exists an $A = (a_{ij})_{1 \leq i, j \leq n} \in \Am(\Bc)(\Z)$ with $A = RU$. Then there exists finitely many choices, only depending on $\Bc$,  of $a_{i, i+1}$ and $a_{i+1, i}$ by the proof of Proposition \ref{IntegralCartan}. Note the choice of the partial column vector
    \[
    \mathbf{v}_i = (a_{i-1, i}, a_{ii}, a_{i+1, i})
    \]
    is finite.
    
    Write $R = (m_{ij})_{1 \leq i \leq 3, 1 \leq j \leq n}$ with
    \[
    \mathbf{m}_j = (m_{1j}, m_{2j}, m_{3j}) = (\Delta_{(2, 3, j)}, -\Delta_{(1, 3, j)}, \Delta_{(1, 2, j)}) \in \Q^3.
    \]
    and its minor
    \[
    M_i = (\mathbf{m}_{i-1}, \mathbf{m}_{i}, \mathbf{m}_{i+1})^t \in \GL(3, \Q).
    \]
    Then $U = (M_i^{-1} \mathbf{v}_i)_{1 \leq i \leq n}$ is uniquely determined by $M$ and chosen $\mathbf{v}_i$. In particular, $a_{ij}$ is given as
    \[
    a_{ij} = \mathbf{m}_i^t (M_j^{-1} \mathbf{v}_j) = \frac{1}{\Delta_{(j-1, j, j+1)}}(\Delta_{(i, j, j+1)}, \Delta_{(j-1, i, j+1)}, \Delta_{(j-1, j, i)})^t \mathbf{v}_j,
    \]
    only depending on the choice of $M$ and $\mathbf{v}_j$.
\end{proof}

\begin{corollary}
    There is a surjection 
    \[
    \Am(\Bc)(\Z) \rightarrow \mathscr{G}_{\Bm} (\Bc)(\Z)
    \]
    with finite fiber. In particular,
    \[
    |\mathcal{C}(\Bm; \Bc)(\Z)| < \infty \Longleftrightarrow |\mathscr{G}_{\Bm} (\Bc)(\Z) / D | < \infty.
    \]
\end{corollary}

\section{Proof of Theorem \ref{mainthm1}}\label{sec:proof-main}
Throughout the section, we assume $n\ge 5$.  Let $\mathcal U_n$ denote the set of matrices $R=(I_3\mid M)$ with $M\in D^+(\mathbf\Delta)$ such that
\begin{align*}
  \Delta_I(R)&\in\mathbb Z_{>0}
  &&\text{for every positive triple }I,\\
  \Delta_i(R)&=1
  &&\text{for }1\leq i\leq n-1.
\end{align*}
The set $\mathcal U_n$ depends only on $n$. For $R=(I_3\mid M)\in\mathcal U_n$, we have
\[
 \Phi_{\mathbf \Delta}(M) = \bigl(\Delta_{I_1}(M),\ldots,\Delta_{I_{2n-8}}(M);1,\ldots,1\bigr).
\]

For $R\in\mathcal U_n$, the identity 
\[
\mathbf{m}_i=(\Delta_{(2, 3, j)}, -\Delta_{(1, 3, j)}, \Delta_{(1, 2, j)})^t
\]
from the proof of Theorem~\ref{IntegralFactorizationII} shows that $R\in \mathcal \Mat_{3\times n}(\Z)$. Moreover, the matrix 
\[
 M_i = (\mathbf{m}_{i-1}, \mathbf{m}_{i}, \mathbf{m}_{i+1})^t \in\SL(3,\Z),
\]
since $\det M_i=\Delta_{i-1, i, i+1}=\Delta_{i-2}=1$.

For each $i\in\Z/n\Z$, put
\[
  \kappa_i=
  \begin{cases}
    4\cos^2(\pi/m_i),&m_i<\infty,\\
    c_i,&m_i=\infty,
  \end{cases}
\]
and assume that every $\kappa_i$ is an integer.

We first fix the diagonal and consecutive entries of the matrix to be reconstructed. Set $a_{ii}=2$, and choose nonpositive integers $a_{i,i+1}$ and $a_{i+1,i}$ satisfying
\begin{equation}\label{integralentries}
  a_{i,i+1}a_{i+1,i}=\kappa_i,\qquad a_{i,i+1}=0\iff a_{i+1,i}=0
\end{equation}
for every $i$. Such a choice always exists: take $(a_{i,i+1},a_{i+1,i})=(0,0)$ if $\kappa_i=0$ and $(-\kappa_i,-1)$ otherwise. These entries are fixed independently of $R$.

For each $i$, let
\[
  v_i=(a_{i-1,i},2,a_{i+1,i})^t.
\]
For $R\in\mathcal U_n$, define the matrix $A(R)$ by
\[
  \bigl(A(R)\bigr)_{ji}:=\mathbf m_j^{\mathsf t}M_i^{-1}v_i \qquad(1\leq j,i\leq n),
\]
as in the reconstruction in the proof of Theorem~\ref{IntegralFactorizationII}.

\begin{proposition}\label{welldefinedness}
    For every $R\in \mathcal U_n$, the reconstructed matrix satisfies
    \[
    A(R)\in\Am(\Bc)(\Z).
    \]
    In particular, $\mathcal U_n\subset \mathscr{G}_{\Bm}(\Bc)(\Z)$.
\end{proposition}
\begin{proof}
    Set $b_i=M_i^{-1}v_i$ and $U=(b_1,\ldots,b_n)$. Then $A(R)=R^{\mathsf t}U$. Then $A=R^tU$, where $U=(b_1,\dots,b_n)$. Since $M_i\in\SL(3,\Z)$ and $v_i\in\Z^3$, each $b_i$ is integral. As $R$ is integral, $A$ has integral entries.
    Moreover, the equations $M_i b_i=v_i$ give the following identities:
    \[
    (A_{i-1,i}(R),A_{ii}(R),A_{i,i+1}(R))^t=(a_{i,i-1},2,a_{i+1,i})^t
    \]
    Thus $A(R)$ agrees with the chosen diagonal and consecutive entries.
    
    Write $A(R)=(a_{ji})_{1\le i,j\le n}$. Since $\Delta_{i-1,i,i+1}=1$, the reconstruction formula in the proof of Theorem~\ref{IntegralFactorizationII} gives
    \[
    a_{ji}=a_{i-1,i}\Delta_{j,i,i+1}+2\Delta_{i-1,j,i+1}+a_{i+1,i}\Delta_{i-1,i,j}.
    \]
    For $j\notin\{i-1,i,i+1\}$, the definition of $\mathcal U_n$ and cyclic order of the indices give
    \[
    \Delta_{j,i,i+1},\Delta_{i-1,i,j}\geq1,\qquad \Delta_{i-1,j,i+1}\leq-1.
    \]
    Since $a_{i-1,i},a_{i+1,i}\le0$, it follows that $a_{ji}\leq-2$.Interchanging $i$ and $j$, we also obtain $a_{ij}\leq-2$. Hence $a_{ij}a_{ji}\geq4$. Moreover formula~\ref{integralentries} implies that $A$ satisfies Cartan matrix conditions.
    
    It remains to realize $A$ as the Cartan matrix of a Coxeter polygon realizing $Q$. Consider the cone
    \[
    C_R=\bigl\{x\in \R^3\mid \mathbf{m}_j^tx\le0\text{ for all }j\bigr\}.
    \]
    For each $i$, put $w_i=-\mathbf{m}_j\times\mathbf{m}_{j+1}$. Then
    \[
    \mathbf{m}_i^tw_i=\mathbf{m}_{i+1}^tw_i=0,\qquad\mathbf{m}_j^tw_i=-\Delta_{j,i,i+1}<0\quad\text{for }j\notin\{i,i+1\}.
    \]
    In particular, $\sum_iw_i$ lies in the interior of $C_R$. Since the columns of $R$ span $\R^3$, $C_R$ contains no line. Thus its projectivization $P_R=\P(C_R\setminus\{0\})$ is a properly convex $n$-gon. For each $i$, the intersection $F_i:=P_R\cap\P(\ker\mathbf{m}_j^t)$ is a side of $P_R$ with endpoints $[w_{i-1}]$ and $[w_i]$. Hence $P_R$ is a $n$-gon with side $F_1,\dots,F_n$ in cyclic order.

    Define $\sigma_i=\mathrm{Id}_3-b_i\mathbf{m}_i^t$. Since $\mathbf{m}_i^tb_i=a_{ii}=2$, each $\sigma_i$ is a reflection fixing $F_i$ pointwise. The Cartan conditions verified above show that $(P,(\sigma_i)_{1\le i\le n})$ realizes $Q$. Thus $A$ is compatible with $Q$. Since $n\ge5$, Fact~\ref{polygon-cartan} gives $A\in\Am$. Finally, the identities $a_{i,i+1}a_{i+1,i}=\kappa_i$ imply that $A\in\Am(\Bc)$.
\end{proof}

\begin{corollary}\label{frozenfinitefibers}
    With the consecutive entries fixed as above, the map
    \[
    \Psi_{\Bm,\Bc}:\mathcal U_n\longrightarrow \Ccal(\Bm;\Bc)(\Z), \quad R\longmapsto[A(R)].
    \]
    is well defined and has finite fibers. 
\end{corollary}
\begin{proof}
    Well-definedness follows from Proposition~\ref{welldefinedness}. For $R\in\mathcal U_n$, write $A=A(R)=R^tU$. Since $\rk A=3$, choose three linearly independent columns of $A$. Let $B=(B_0\mid B_1)$ be the transpose of the resulting $n\times 3$ matrix. Since $R=(I_3\mid M)$, the factorization gives $B=B_0R$. The condition $\rk B=3$ implies that  $B_0$ is invertible, and hence $R=B^{-1}_0B$. Thus $A$ determines $R$ uniquely. 
    
    By Proposition~\ref{IntegralCartan}, each positive diagonal conjugacy class contains only finitely many Cartan matrices with integer entries. This shows that $\Psi_{\Bm,\Bc}$ has finite fibers.
\end{proof}

\begin{proposition}\label{infiniteU_n}
    The set $\mathcal U_n$ is infinite for every $n\ge9$.
\end{proposition}
\begin{proof}
    By Scott \cite{scott2006grassmannians}, the cluster algebra $\Acal_{3,n}$ has finitely many cluster variables if and only if $n < 9$. Suppose $n \geq 9$ and fix an initial cluster $\xb = (x_1, \ldots, x_{2n-8}; \Delta_0, \ldots, \Delta_{n-1})$. By the Laurent phenomenon, the image of any cluster $\xb'$ by the evaluation $\xb \mapsto (1, \ldots, 1; 1, \ldots, 1)$ is an integral vector corresponding to an $M \in \mathcal U_n$.

    It suffices to show that the map always has a finite fiber. Given an integral vector, by the positivity phenomenon, there exists finitely many Laurent polynomial tuples. Thus the finiteness of $\mathcal U_n$ is equivalent to the finiteness of cluster variables of $\Acal_{3,n}$.
\end{proof}

\begin{theorem}\label{infinite Cartan}
    Suppose that $n\ge 9$ and that $\kappa_i\in\Z$ for every $i$. Then $\Ccal(\Bm;\Bc)(\Z)$ is infinite.
\end{theorem}
\begin{proof}
    By Corollary~\ref{frozenfinitefibers}, the map $\Psi_{\Bm,\Bc}\colon\mathcal U_n\to\Ccal(\Bm;\Bc)(\Z)$ has finite fibers. Since $\mathcal U_n$ is infinite by Proposition~\ref{infiniteU_n}, its image is infinite, and hence so is $\Ccal(\Bm;\Bc)(\Z)$.
\end{proof}

\section{Analysis for small polygons}\label{sec:proof-main2}
\subsection{Finiteness for pentagons}
Let $Q$ be a labeled pentagon. Consider a Cartan matrix $A$ compatible with a labeled pentagon $Q$, with factorization $A=R^tU$. Since $\rk A=\rk R=\rk U=3$, we have
\[
\ker A^t=\ker R,\qquad \ker A=\ker U.
\]
Consider the basis matrices $B_1$ of $\ker R$ and $B_2$ of $\ker U$, respectively.  For $i<j$, let $\Delta_{ij}(B_\ell)$ denote the corresponding two-row minor. By the complementary minor identity, there exist nonzero $c_R,c_U\in \R$ such that for $1\le i<j\le5$ and $J=\{1,2,3,4,5\}\setminus\{i,j\}$,
\[
\Delta_{ij}(B_1)=c_R(-1)^{i+j}\Delta_{J}(R),\qquad\Delta_{ij}(B_2)=c_U(-1)^{i+j}\Delta_{J}(U).
\]

The proper convexity of $P_A$ implies that, after a suitable choice of orientation, $\Delta_I(R)>0$ for every three-element set $I$. Since each such $I$ contains a nonadjacent pair $i,j$, we have
\[
\Delta_I(R)\Delta_I(U)=\det A_{I,I}\le 8-2a_{ij}a_{ji}<0,
\]
so $\Delta_I(U)<0$. Setting $D=\diag(1,-1,1,-1,1)$, the complementary-minor identities show that all ordered $2\times2$ minors of each of $DB_1$ and $DB_2$ have the same nonzero sign. By choosing the two kernel bases appropriately, we may assume that all these minors are positive.

Let $A,A'$ be compatible Cartan matrices of a labeled pentagon with realization $A=R^tU$ and $A'=(R')^tU'$. Let $B_1$ and $B_2'$ be as above. The positivity of the ordered minors of $DB_1$ and $D'B_2'$ determines the cyclic order of their row directions in $\mathbb{RP}^1$. Then we may write
\[
B_1=\Lambda K(t,T),\qquad B_2'=WK(s,S),\qquad K(t,T)=\begin{pmatrix}
    1&0\\
    -1&-t\\
    1&1\\
    -1&-T\\
    0&1
\end{pmatrix}
\]
where $0<t,s<1$, $1<T,S$, and $W=\diag(w_i)$ and $\Lambda=\diag(\lambda_i)$ are positive diagonal matrices.

Set $\beta=t+\frac{1-t}{S}$, $\gamma=s+\frac{1-s}{T}$ and define
\[X_1=B_1,\qquad Y_2=B_2'\begin{pmatrix}
    1&-\gamma\\
    -\beta&1
\end{pmatrix}.\] Since $t<\beta<1$ and $s<\gamma<1$, $X_1$ and $Y_2$ are still basis matrices of $\ker R$ and $\ker U'$, respectively.

Writing their rows as $x_i^t$ and $y_i^t$, direct computation gives
\begin{equation}\label{inner_product_positive}
\langle x_i,y_j\rangle>0
\qquad\text{whenever $i,j$ are nonadjacent}.
\end{equation}
We shall use this choice of bases in the proof below.

\begin{theorem}
    For $\Bm=(m_1,m_2,m_3,m_4,m_5)$ and $\Bc=(c_i)_{i\in I_\infty}$, the set $\Ccal(\Bm;\Bc)(\Z)$ is finite.
\end{theorem}
\begin{proof}
    Since the adjacent products are fixed by $\Bm$ and $\Bc$, there are only finitely many choices of integral adjacent entries. Fix one such choice.

    Suppose $A,A'\in\Ccal(\Bm;\Bc)(\Z)$ with realization $A=R^tU$, $A'=(R')^tU'$. Suppose that $A\le A'$ entrywise and put $E=A'-A\ge0$. Then $E$ is supported on nonadjacent entries. Choose basis matrices $X_1$ of $\ker R$ and $Y_2$ of $\ker U'$ satisfying \eqref{inner_product_positive}.
    We have
    \[
    X_1^tEY_2=X_1^t(A'-A)Y_2=0.
    \]
    Taking the trace gives
    \[
    0=\tr(X_1^tEY_2)
    =\sum_{i,j}E_{ij}\langle x_i,y_j\rangle.
    \]
    The summand is nonnegative implies that $E=0$ and $A=A'$.

    Consequently, for the fixed adjacent entries, the vectors $(-a_{ij})_{i-j\equiv\pm2\;(\mathrm{mod}\;5)}$ form an antichain in $\Z^{10}_{>0}$. By Dickson's Lemma(\cite[Lemma~A]{Dickson1913}), the set of tuple $(a_{ij})$ is finite in $\Z^{10}_{>0}$. Since there are only finitely many choices of adjacent entries, the claim follows from Propositions~\ref{Integrality criterion} and~\ref{IntegralCartan}.
\end{proof}

\subsection{An infinite family for symmetric hexagon}

Consider the following forms of Cartan matrices
\[
 A(a,b,c,d)=
 \begin{pmatrix}
 2&-1&-a&-c&-b&-1\\
 -1&2&-1&-b&-d&-a\\
 -b&-1&2&-1&-a&-c\\
 -d&-a&-1&2&-1&-b\\
 -a&-c&-b&-1&2&-1\\
 -1&-b&-d&-a&-1&2
 \end{pmatrix}.
\]

Order the sides as $(1,3,5,4,6,2)$, write $J_3$ for the all-ones
$3\times3$ matrix, and set
\[
 X=\begin{pmatrix}2&-a&-b\\-b&2&-a\\-a&-b&2\end{pmatrix},\qquad
 C=(c-1)\Id+J_3,\qquad D=(d-1)\Id+J_3.
\]
Observe $C, D$ commutes with $X$. Suppose that $D$ invertible\footnote{This holds if and only if $d \neq -2, 1$.} and $XX^t = CD$ then $XD^{-1}X^t = C$ and
\[
 \widetilde A=\begin{pmatrix}X&-C\\-D&X^t\end{pmatrix}=
 \begin{pmatrix}X\\-D\end{pmatrix}
 \begin{pmatrix}\Id&-D^{-1}X^t\end{pmatrix},
 \qquad \rank\widetilde A=3.
\label{eq:hexfactor}
\]

Since $XX^t = (a^2 + b^2 - ab + 2a + 2b + 4)\Id + (ab - 2a - 2b)J_3$, it suffices to solve a Diophantine system
\[
\begin{cases}
 (c-1)(d-1) = a^2 + b^2 - ab + 2a + 2b + 4 \\
 c+d+1 = ab - 2a - 2b
\end{cases}
\]

Consequently we have a Diophantine system
\begin{equation}\label{HexDiophantine}
    c+d = ab - 2a - 2b - 1, \quad cd = a^2 + b^2 +2
\end{equation}

\begin{theorem}\label{HexInfinite}
For $\mathbf{m}=(3,3,3,3,3,3)$, the set $\mathcal{C}(\mathbf{m})(\Z)$ is infinite.
\end{theorem}

\begin{proof}
    It suffices to check that the Diophantine system \eqref{HexDiophantine} has infinitely many positive integral solutions with $d \neq 1$. Put $x = a-2$, $y = b-2$ and $z = c$, then we obtain
    \[
    d= xy - z - 5
    \]
    by the first equation. Then the second equation is written as
    \[
    x^2 + y^2 + z^2 -xyz + 4x + 4y + 5z + 10 = 0
    \]
    which is a generalized Markoff equation, having three Vieta involutions
    \begin{align*}
    V_x : (x, y, z) &\rightarrow (yz - 4 - x, y, z) \\
    V_y : (x, y, z) &\rightarrow (x, xz - 4 - y, z) \\
    V_z : (x, y, z) &\rightarrow (x, y, xy - 5 - z).
    \end{align*}
    Write $\Gamma = \langle V_x, V_y, V_z \rangle$. We can determine all (nondegenerate) fundamental solutions of the equation and whether it has an infinite $\Gamma$-orbit or not \cite{whang2020nonlinear, biswas2022surface}. There exists a fundamental solution
    \[
    (x, y, z) = (6, 6, 5)
    \]
    whose coordinates are positive integers greater than $2$. By induction one can check that the $\Gamma$-orbit of each fundamental solution gives a tuple of positive integers greater than $2$. 
\end{proof}

\begin{remark}
    For any \textit{symmetric} labelling $\mathbf{m} = (k, k, k, k, k, k)$ (if $k = \infty$, we need all of the $c_i$s are equal) unless $k = 2$, the proof goes similarly and one can determine whether the generalized Markoff equation has a fundamental solution or not, in finite time.
\end{remark}

\subsection{An infinite family for symmetric octagon}

Consider the following forms of Cartan matrices
\[
 A(a,b,c,d,e)=
 \begin{pmatrix}
 2&-1&-a&-c&-e&-c&-b&-1\\
 -1&2&-1&-b&-d&-e&-d&-a\\
 -b&-1&2&-1&-a&-c&-e&-c\\
 -d&-a&-1&2&-1&-b&-d&-e\\
 -e&-c&-b&-1&2&-1&-a&-c\\
 -d&-e&-d&-a&-1&2&-1&-b\\
 -a&-c&-e&-c&-b&-1&2&-1\\
 -1&-b&-d&-e&-d&-a&-1&2
 \end{pmatrix}.
\]

Order the sides as $(1,3,5,7,2,4,6,8)$, write
\[
 P = 
 \begin{pmatrix}
     0&1&0&0\\
     0&0&1&0\\
     0&0&0&1\\
     1&0&0&0
 \end{pmatrix}
\]
then set
\[
        X = 2\Id - aP - eP^2 - bP^3, \quad
    C = \Id + cP + cP^2 + P^3, \quad
    D = \Id + P + dP^2 + dP^3.
\]
The upshot is to construct a Diophantine equation $XX^t = CD$. Observe $X, C, D$ commutes each other and $C, D$ has a common kernel generated by $\mathbf v = (1, -1, 1, -1)$, the eigenvector of $P$ of eigenvalue $-1$. Put $e = a+b+2$, then $X, C, D$ has the simultaneous kernel generated by $\mathbf v$. Write $W = \langle \mathbf v \rangle ^\perp$ and suppose that $D|_W$ is invertible\footnote{This holds if and only if $d \neq \pm 1$.} then
\[
 \widetilde A|_W=\begin{pmatrix}X|_W&-C|_W\\-D|_W&(X|_W)^t\end{pmatrix}=
 \begin{pmatrix}X|_W\\-D|_W\end{pmatrix}
 \begin{pmatrix}\Id&-D|_W^{-1}(X|_W)^t\end{pmatrix},
 \qquad \rank\widetilde A=3.
\label{eq:hexfactor}
\]
if and only if $XX^t = CD$. Consequently we have a Diophantine system
\begin{equation}
    x^2 + y^2 + z^2 - xyz + 8x + 8y + 10z + 32 = 0
\end{equation}
by letting $x = a-2$, $y = b-2$, $z = c-1$ as before. By taking a fundamental solution $(8, 8, 6)$ and similarly use the method of the proof of Theorem \ref{HexInfinite}, we have the following theorem.

\begin{theorem}
    For $\mathbf{m}=(3,3,3,3,3,3,3,3)$, the set $\mathcal{C}(\mathbf{m})(\Z)$ is infinite.
\end{theorem}


\printbibliography

\end{document}